\documentclass[secthm,dvips]{arxbj}

\aid{0}
\volume{16}
\issue{3}
\pubyear{2010}
\firstpage{759}
\lastpage{779}
\doi{10.3150/09-BEJ232}

\makeatletter
\newcommand{\argmax}{\operatorname{arg\,max}}
\newtheorem{lema}{Lemma}[section]
\newtheorem{assump}{Assumption}[section]
\makeatother

\begin{document}
\begin{frontmatter}

\title{Asymptotic properties of maximum likelihood estimators in models
with multiple change points}
\runtitle{Multiple-change-point models}

\begin{aug}
\author[1]{\fnms{Heping} \snm{He}\thanksref{1}\ead[label=e1]{hhe@math.ku.edu}}
\and
\author[2]{\fnms{Thomas A.} \snm{Severini}\thanksref{2}\ead[label=e2]{severini@northwestern.edu}}
\runauthor{H. He and T.A. Severini}
\address[1]{Department of Mathematics, University of Kansas, 1460
Jayhawk Blvd, Lawrence,
KS 66045, USA. \printead{e1}}
\address[2]{Department of Statistics, Northwestern University,
Evanston, IL 60208, USA.\\ \printead{e2}}
\end{aug}

\received{\smonth{5} \syear{2008}}
\revised{\smonth{8} \syear{2009}}

%
\begin{abstract}
Models with multiple change points are used in many fields; however,
the theoretical properties of maximum likelihood estimators of such
models have received relatively little attention. The goal of this
paper is to establish the asymptotic properties of maximum likelihood
estimators of the parameters of a multiple change-point model for a
general class of models in which the form of the distribution can
change from segment to segment and in which, possibly, there are
parameters that are common to all segments. Consistency of the maximum
likelihood estimators of the change points is established and the rate
of convergence is determined; the asymptotic distribution of the
maximum likelihood estimators of the parameters of the within-segment
distributions is also derived.
Since the approach used in single change-point models
is not easily extended to multiple change-point models,
these results require the introduction of those tools for analyzing the
likelihood function in a multiple change-point model.
\end{abstract}

%
\begin{keyword}
\kwd{change-point fraction}
\kwd{common parameter}
\kwd{consistency}
\kwd{convergence rate}
\kwd{Kullback--Leibler distance}
\kwd{within-segment parameter}
\end{keyword}

\end{frontmatter}
%
\section{Introduction}\label{sec1}
A change-point model for a sequence of
independent random variables $X_1, \ldots, X_n$
is a model in which there exist unknown change points $n_1, \ldots,
n_k$, $0 = n_0 < n_1 < \cdots< n_k < n_{k+1} = n$, such that, for
each $j=1,2,\ldots, k+1$, $X_{n_{j-1}+1}, \ldots, X_{n_j}$\vspace*{1pt} are
identically distributed with a distribution that depends on $j$.
Here, we consider parametric change-point models in which the
distribution of $X_{n_{j-1} +1}, \ldots, X_{n_j}$ is parametric;
however, the form of the distribution can be different for each $j$.
Change-point models are used in many fields. For example,
Broemeling and Tsurumi~(\citeyear{Broemeling87}) uses a multiple change-point model for
the US demand
for money; Lombard~(\citeyear{Lombard86}) uses a multiple change-point model to model
the effect of sudden changes in wind direction on the flight of a
projectile; Reed (\citeyear{Reed98}) uses a multiple change-point model in the
analysis of forest fire data. A~number of authors have used multiple
change-point models in the analysis of DNA sequences; see, for
example, Braun and Muller (\citeyear{Braun98}), Fu and Curnow (\citeyear{Fu90a,Fu90b})
and Halpern~(\citeyear{Halpern00}). Many further examples are provided in
the monographs Chen and Gupta (\citeyear{Chen00}) and Cs\"{o}rg\H{o} and
Horv\'{a}th (\citeyear{Csorgo97}).

The goal of this paper is to establish the asymptotic properties of
maximum likelihood estimators of the parameters of a multiple
change-point model, under easily verifiable conditions.
These results are based on the following model.
Assume that the vectors in the data set $x_1,x_2,\ldots,x_n$ are
independently drawn from the parametric model
\[
f_j(\psi^0,\theta_j^0;x_i),\qquad  n_{j-1}^0+1\leq i\leq
n_j^0,  j=1,2,\ldots,k+1,
\]
where $f_j(\psi^0,\theta_j^0;x)$ is a probability density function of a
continuous distribution with unknown common parameter $\psi^0$ for all
$j=1,2,\ldots,k+1$ and unknown within-segment parameters $\theta_j^0$
for each $j=1,2,\ldots,k+1$; $f_j(\psi^0,\theta_j^0;x)$ may have the
same functional form for some or all of $j=1,2,\ldots,k+1$; $\psi^0$ may
be a vector; $\theta_j^0$ may be a different vector parameter of
different dimensions for each $j=1,2,\ldots,k+1$. In this model, there
are $k$ unknown change points~$n_1^0,n_2^0,\ldots,n_k^0$, where the
number of change points $k$ is assumed to be known. The parameter $\psi
^0$ is common to all segments.

There are a number of results available on the asymptotic properties of
parameter estimators in change-point
models. See, for example, Hinkley~(\citeyear{Hinkley70,Hinkley72}), Hinkley and
Hinkley~(\citeyear{Hinkleyetal70}), Battacharya (\citeyear{Battacharya87}), Fu and Curnow (\citeyear{Fu90a,Fu90b}),
Jandhyala and Fotopoulos (\citeyear{Jandhyala99,Jandhyala01}) and Hawkins (\citeyear{Hawkins01}); the two monographs Chen and Gupta
(\citeyear{Chen00}) and Cs\"{o}rgo and Horv\'{a}th (\citeyear{Csorgo97})
have detailed bibliographies on this topic.

In particular, Hinkley (\citeyear{Hinkley70}) considers likelihood-based inference for
a single change-point model, obtaining the asymptotic distribution of
the maximum likelihood estimator of the change point under the
assumption that the other parameters in the model are known. Hinkley
(\citeyear{Hinkley70}) and Hinkley (\citeyear{Hinkley72}) argue
that this asymptotic distribution is also valid when the parameters are unknown.

Unfortunately, there are problems in extending the approach used in
Hinkley (\citeyear{Hinkley70,Hinkley72}) to the setting considered here.
The method used in Hinkley (\citeyear{Hinkley70,Hinkley72}) is based on
considering the relative locations of a candidate change point and the
true change point. When there is only a single change point, there are
only three possibilities: the candidate change point is either
greater than, less than or equal to the true change point. However, in
models with $k$ change points, the relative
positions of the candidate change points and the true change points can
become quite complicated and the
simplicity and elegance of the single change point argument is lost.

A second problem arises when extending the
argument for the case in which the change points are the only
parameters in the model to the case in which there are unknown
within-segment parameters. The consistency argument used in the former
case is extended to the latter case
using a ``consistency assumption'' (Hinkley (\citeyear{Hinkley72}), Section 4.1); this
condition is discussed in Appendix \hyperref[appA]{A} and examples
are given which show that this assumption is a strong one that is not
generally satisfied in the class of models considered here.

There are relatively few results available on the asymptotic properties
of maximum likelihood estimators in multiple
change-point models. Thus, the present paper has done several things.
In the general model described above,
in which there is a fixed, but arbitrary, number of change points, we
show that the maximum likelihood estimators
of the change points are consistent and converge to the true change
points at the rate $1/n$, under relatively weak regularity conditions.
As noted above,
a simple extension of the approach used in single change-point models
is not available; thus, the second thing achieved by this paper is the
introduction of the tools necessary for analyzing the likelihood
function in a multiple change-point
model. Finally, the asymptotic distribution of the maximum likelihood
estimators of the parameters of the
within-segment distributions is derived for the general case described
above, in which the form of the
distribution can change from segment to segment and in which, possibly,
there are parameters that are
common to all segments.

The paper is organized as follows. The asymptotic theory of maximum
likelihood estimators of a multiple change-point model is described in
Section \ref{sec2}. Section \ref{sec3} contains a numerical example illustrating these results
and Section \ref{sec4} contains some discussion of future research which builds
on the results given in this paper.
Appendix \hyperref[appA]{A} discusses the ``consistency assumption'' used in Hinkley
(\citeyear{Hinkley72}); all technical proofs are given in Appendix \hyperref[appB]{B}.

\section{Asymptotic theory}\label{sec2}
Consider estimation of the multiple
change-point model introduced in Section \ref{sec1}. For any change point
configuration $0=n_0<n_1<n_2<\cdots<n_{k}<n_{k+1}=n$, the
log-likelihood function is given by
\[
l\equiv l(n_1, \ldots, n_k, \theta_1, \ldots, \theta_{k+1}, \psi) =\sum
_{j=1}^{k+1}\sum_{i=n_{j-1}+1}^{n_j}\log f_j(\psi,\theta_j;x_i).
\]
Estimators of all change points, all within-segment parameters and the
common parameter are given by
\[
(\hat{n}_1,\hat{n}_2,\ldots,\hat{n}_k,\hat{\theta}_1,\hat{\theta}_2,\ldots
,\hat{\theta}_{k+1},\hat{\psi})=\mathop{\argmax}_{0<n_1<n_2<\cdots<n_k<n; \theta_j\in\Theta_j,j=1,2,\ldots
,k+1;\psi\in\Psi} l,
\]
where $\Theta_j, j=1,2,\ldots,k+1$, and $\Psi$ are the parameter spaces
of $\theta_j$, $j=1, \ldots, k+1$, and $\psi$, respectively.

Let
\begin{eqnarray*}
\lambda_j^0&=&n_j^0/n \qquad\mbox{for } j=1,2,\ldots
,k,\\
\lambda_j&=&n_j/n \qquad\mbox{for } j=1,2,\ldots,k
,\\
\lambda^0&=&(\lambda_1^0,\lambda_2^0,\ldots,\lambda_k^0),\\
\lambda&=&(\lambda_1,\lambda_2,\ldots,\lambda_k);\nonumber\\
\theta^0&=&(\theta_1^0,\theta_2^0,\ldots,\theta_{k+1}^0),\\
\theta&=&(\theta_1,\theta_2,\ldots,\theta_{k+1}) ,\\
\phi^0&=&(\psi^0,\theta^0)=(\psi^0,\theta_1^0,\theta_2^0,\ldots,\theta
_{k+1}^0) ,\\
\phi&=&(\psi,\theta)=(\psi,\theta_1,\theta_2,\ldots,\theta
_{k+1}).
\end{eqnarray*}
Note that $\lambda^0$ is taken to be a constant vector as $n$ goes to infinity.

Define
\vspace*{-6pt}
\begin{eqnarray*}
\hat{\ell}^{(j)}(\psi,\theta_j)&=&\sum_{i=\hat{n}_{j-1}+1}^{\hat
{n}_j}\log f_j(\psi,\theta_j;x_i),\qquad j=1,2,\ldots,k+1,
\\
\ell^{(j)}(\psi,\theta_j)&=&\sum_{i=n^0_{j-1}+1}^{n^0_j}\log f_j(\psi
,\theta_j;x_i),\qquad j=1,2,\ldots,k+1,
\\
\hat{\ell}(\psi,\theta)&=&\sum_{j=1}^{k+1}\sum_{i=\hat
{n}_{j-1}+1}^{\hat{n}_j}\log f_j(\psi,\theta_j;x_i) ,
\\
\ell^0(\psi,\theta)&=&\sum_{j=1}^{k+1}\sum
_{i=n^0_{j-1}+1}^{n^0_j}\log f_j(\psi,\theta_j;x_i) ,
\\
\ell(\psi,\theta)&=&\sum_{j=1}^{k+1}\sum
_{i=n_{j-1}+1}^{n_j}\log f_j(\psi,\theta_j;x_i) .
\end{eqnarray*}

The expected information matrix is given by
\[
i(\psi,\theta)=E[-\ell^0_{\phi\phi}(\psi,\theta);\phi]
=\pmatrix{
E[-\ell^0_{\psi\psi}(\psi,\theta);\phi]& E[-\ell^0_{\psi\theta}(\psi
,\theta);\phi]\vspace*{3pt}\cr
E[-\ell^0_{\psi\theta}(\psi,\theta);\phi]^T&E[-\ell^0_{\theta\theta
}(\psi,\theta);\phi]},
\]
\vspace*{-6pt}
\begin{eqnarray*}
&&E[-\ell^0_{\psi\theta}(\psi,\theta);\phi]\\
&&\quad=\bigl(E\bigl[-\ell^{(1)}_{\psi\theta_1}(\psi,\theta_1);\phi\bigr],E\bigl[-\ell
^{(2)}_{\psi\theta_2}(\psi,\theta_2);\phi\bigr],\ldots,E\bigl[-\ell^{(k+1)}_{\psi
\theta_{k+1}}(\psi,\theta_{k+1});\phi\bigr]\bigr),
\\
&&E[-\ell^0_{\theta\theta}(\psi,\theta);\phi]\\
&&\quad =\operatorname{diag}\bigl(
E\bigl[-\ell^{(1)}_{\theta_1\theta_1}(\psi,\theta_1);\phi\bigr],
E\bigl[-\ell^{(2)}_{\theta_2\theta_2}(\psi,\theta_2);\phi\bigr], \ldots,
E\bigl[-\ell^{(k+1)}_{\theta_{k+1}\theta_{k+1}}(\psi,\theta_{k+1});\phi\bigr]\bigr),
\end{eqnarray*}
where $\operatorname{diag}(\cdot)$ denotes a diagonal block matrix whose diagonal
blocks are in the bracket, other elements are zeros
and the average expected information matrix is given by
\[
\bar{i}(\psi,\theta)=\lim_{n\rightarrow\infty}\frac{1}{n}i(\psi,\theta
) .
\]

The asymptotic properties of these estimators are based on the
following regularity conditions. Other than the parts concerning change
points, these conditions are typically similar to those required for
the consistency and asymptotic normality of maximum likelihood
estimators of parameters in models without change points; see, for
example, Wald (\citeyear{Wald49}). Particularly, compactness of parameter spaces is
a common assumption in the classical likelihood literature.

These conditions are different from those required by Ferger (\citeyear{Ferger01}) and
D\"{o}ring (\citeyear{Doring07}), who consider estimation of change points
in a nonparametric setting in which nothing is assumed about the
within-segment distributions, using a type of nonparametric M-estimator
based on empirical processes.
Thus, these authors do not require conditions on the within-segment
likelihood functions; on the other hand, their method does not provide
estimators of
within-segment parameters.\vspace*{-1pt}

\begin{assump}\label{a20}
It is assumed that
for $j=1,2,\ldots,k$, $f_{j+1}(\psi^0,\theta_{j+1}^0;x)\neq
f_{j}(\psi^0,\theta_{j}^0;x)$ on a~set of non-zero measure.
\end{assump}

This assumption guarantees that the distributions in two neighboring
segments are different; clearly, this
is required for the change points to be well defined.\vspace*{-1pt}

\begin{assump}\label{a21}
It is assumed that:
\begin{longlist}[1.]
\item[1.] for $j=1,2,\ldots,k+1$, $\theta_j$ and $\theta_j^0$ are contained in
$\Theta_j$, where $\Theta_j$ is a compact subset of~$\mathcal
{R}^{d_j}$; $\psi$ and $\psi^0$ are contained in $\Psi$ where $\Psi$ is
a compact subset of $\mathcal{R}^{d}$; here, $d, d_1, \ldots, d_{k+1}$
are non-negative integers;
\item[2.] $\ell(\psi,\theta)$ is third-order continuously differentiable
with respect to $\psi,\theta$;
\item[3.] the expectations of the first and second order derivatives of $\ell
^0(\psi,\theta)$ with respect to $\phi$ exist for $\phi$ in its
parameter space.
\end{longlist}
\end{assump}

Compactness of the parameter space is used to establish the consistency
of the maximum likelihood estimators
of $n_1/n, \ldots, n_k/n, \theta_1, \ldots, \theta_{k+1}, \psi$; see, for
example, Bahadur (\citeyear{Bahadur71}) for further discussion of this condition and
its necessity in general models. If we assume further conditions on
models, the compactness of the parameter space may be avoided. But this
appears to be a substantial task for future work.
Differentiability of the log-likelihood function is used to justify
certain Taylor series expansions. Both parts of Assumption \ref{a21} are
relatively weak and are essentially the same as conditions used in
parametric models without change points; see, for example, Schervish
(\citeyear{Schervish95}), Section 7.3. Part 3 is very weak and is used in the proof of
Theorem~\ref{thm3}.\vspace*{-1pt}

\begin{assump}\label{a22}It is assumed that:
\begin{longlist}[1.]
\item[1.] for any $j=1,2,\ldots,k+1$ and any integers $s,t$ satisfying $0\leq
s<t\leq n$,
\[
E\Biggl\{\max_{\psi\in\Psi,\theta_j\in\Theta_j}\Biggl(\sum_{i=s+1}^{t}\{\log
f_j(\psi,\theta_j;X_i)-E[\log f_j(\psi,\theta_j;X_i)]\}\Biggr)^2\Biggr\}\leq
C(t-s)^r ,
\]
where $r<2$ and $C$ is a constant;
\item[2.] for any $j=1,2,\ldots,k+1$ and any integers $s,t$ satisfying
$n_{j-1}^0\leq s<t\leq n_j^0$,
\begin{eqnarray*}
&&E\Biggl\{\max_{\psi\in\Psi,\theta_j\in\Theta_j}\Biggl(\sum_{i=s+1}^t\{[\log
f_j(\psi,\theta_j;X_i)-\log f_j(\psi^0,\theta^0_j;X_i)]
-v(\psi,\theta_j;\psi^0,\theta^0_j)\}\Biggr)^2\Biggr\}\\
&&\quad\leq D(t-s)^r ,
\end{eqnarray*}
where $v(\psi,\theta_j;\psi^0,\theta^0_j)$ is introduced in equation
\textup{(\ref{eqn2.2})}, $r<2$ and $D$ is a constant.
\end{longlist}
\end{assump}

Parts 1 and 2 of Assumption \ref{a22} are technical requirements on the
behavior of the log-likelihood function between and within segments,
respectively.
This condition is used to ensure that the information regarding the
within- and between-segment parameters grows quickly enough to
establish consistency and asymptotic normality of the parameter
estimators. These conditions
are relatively weak; it is easy to check that they are satisfied by at
least all distributions in the exponential family.
Consider a probability density function of exponential family form:
\[
f(\eta,x)=h(x)c(\eta)\exp\Biggl(\sum_{i=1}^{m}w_{i}(\eta)t_{i}(x)\Biggr).
\]
It is then straightforward that the Schwarz inequality gives
\begin{eqnarray*}
&&\Biggl(\sum_{i=s+1}^{t}\{\log f(\eta,X_i)-E[\log f(\eta,X_i)]\}\Biggr)^2
\\
&&\quad \leq\Biggl[1+\sum_{q=1}^{m}w_{q}(\eta)^2\Biggr]\\
&&\qquad{}\times\Biggl\{
\Biggl[\sum_{i=s+1}^{t}\bigl(\log h(X_i)-E(\log h(X_i))\bigr)\Biggr]^2+\sum_{q=1}^{m}\Biggl[\sum
_{i=s+1}^{t}\bigl(t_{q}(X_i)-E(t_{q}(X_i))\bigr)\Biggr]^2\Biggr\}.
\end{eqnarray*}
Therefore, Part 1 of Assumption \ref{a22} is satisfied with $r=1$ because the
function $w_{q}(\eta)$ assumed to be continuous can achieve its maximum
on the compact parameter space. Similarly, Part 2 of Assumption \ref{a22}
is also satisfied with $r=1$.

The main results of this paper are given in the following three theorems.
\begin{thm}[(Consistency)]\label{thm1}
Under Assumption \textup{\ref{a20}}, Part 1 of Assumption \textup{\ref{a21}}\vspace*{1pt} and Part 1 of
Assumption \textup{\ref{a22}},
$\hat{\lambda}_i\rightarrow_p\lambda_i^0, \hat{\theta}_j\rightarrow
_p\theta_j^0$
and $\hat{\psi}\rightarrow_p\psi^0$ as $n\rightarrow+\infty$, that is,
$\hat{\lambda}_i-\lambda_i^0=\mathrm{o}_p(1), \hat{\theta}_j-\theta_j^0=\mathrm{o}_p(1)$ and
$\hat{\psi}-\psi^0=\mathrm{o}_p(1)$, where $\hat{\lambda}_i=\hat{n}_i/n$ for
$i=1,2,\ldots,k$ and $j=1,2,\ldots,k+1$.
\end{thm}

Note that $\hat{n}_i, i=1, 2, \ldots, k$, are not consistent (Hinkley
(\citeyear{Hinkley70})); it is the estimators of the change-point fractions $\hat
{\lambda}_i, i=1, 2, \ldots, k$, that are consistent. The consistency of
$\hat{\theta}_j, j=1, 2, \ldots, k+1$, and $\hat{\psi}$ is the same as
the corresponding result in classical likelihood theory for
independent, identically distributed data.

\begin{thm}[(Convergence rate)]\label{thm2}
Under Assumptions \textup{\ref{a20}--\ref{a22}}, we have
\begin{eqnarray}
\lim_{\delta\rightarrow\infty}\lim_{n\rightarrow\infty}P_r(n\|\hat
{\lambda}-\lambda^0\|_{\infty}\geq\delta)=0,\nonumber
\end{eqnarray}
where $\hat{\lambda}=(\hat{\lambda}_1,\hat{\lambda}_2,\ldots,\hat{\lambda}_k),$
$\|\hat{\lambda}-\lambda^0\|_{\infty}=\max_{1\leq j\leq k}|\hat{\lambda
}_j-\lambda^0_j|.$ That is, $\hat{\lambda}_i-\lambda^0_i=\mathrm{O}_p(n^{-1})$
for $i=1, 2, \ldots, k$.
\end{thm}

We now consider the asymptotic distribution of $\hat{\phi}$, where $\hat
{\phi}=(\hat{\psi},\hat{\theta})$.

\begin{thm}[(Limiting distributions)]\label{thm3}
Under Assumptions \textup{\ref{a20}--\ref{a22}},
\begin{eqnarray*}
&&\sqrt{n}(\hat{\phi}-\phi^0)\stackrel{\mathcal{D}}{\longrightarrow}
N_{d+d_1+d_2+\cdots+d_{k+1}}(0,\bar{i}(\psi^0,\theta^0)^{-1}),
\end{eqnarray*}
where $N_{d+d_1+d_2+\cdots+d_{k+1}}(0,\bar{i}(\psi^0,\theta^0)^{-1})$
is the $d+d_1+d_2+\cdots+d_{k+1}$-dimensional multivariate normal
distribution with mean vector zero and covariance matrix $\bar{i}(\psi
^0,\theta^0)^{-1}$.
\end{thm}

The proofs of Theorems \ref{thm1}--\ref{thm3} are based on the following approach.

Define a function $J$ by
%
\begin{eqnarray}\label{eqn2.1}
J&=&\sum_{j=1}^{k+1}\sum_{i=1}^{k+1}\frac{n_{ji}}{n}\biggl\{\int_{-\infty
}^{+\infty}[\log f_j(\psi,\theta_j;x)-\log f_i(\psi^0,\theta
_i^0;x)]f_i(\psi^0,\theta_i^0;x)\,dx\biggr\}\nonumber\\
&&{}+\frac{1}{n}\sum_{j=1}^{k+1}\sum_{i=n_{j-1}+1}^{n_j}\{
\log f_j(\psi,\theta_j;x_i)-E[\log f_j(\psi,\theta_j;X_i)]\}\\
&&{}-\frac{1}{n}\sum_{j=1}^{k+1}\sum_{i=n_{j-1}^0+1}^{n_j^0}\{
\log f_j(\psi^0,\theta_j^0;x_i)-E[\log f_j(\psi^0,\theta_j^0;X_i)]\}
,\nonumber
\end{eqnarray}
where $n_{ji}$ is the number of observations in the set
$[n_{j-1}+1,n_j]\cap[n_{i-1}^0+1,n_i^0]$ for $i,j=1,2,\ldots,k+1$.
We obviously have that
\[
\mathop{\argmax}_{0<n_1<n_2<\cdots<n_k<n; \theta
_j\in\Theta_j,1\leq j\leq k+1;\psi\in\Psi}l=\mathop{\argmax}_{0<n_1<n_2<\cdots<n_k<n;\theta_j\in\Theta_j,1\leq j\leq
k+1;\psi\in\Psi}J;
\]
thus, the maximum likelihood estimators may be defined as the
maximizers of $J$ rather than as the
maximizers of $l$.\vadjust{\goodbreak}

Let $v(\psi,\theta_j;\psi^0,\theta_i^0)$ be defined by
%
\begin{eqnarray}\label{eqn2.2}
&&v(\psi,\theta_j;\psi^0,\theta_i^0)=\int_{-\infty}^{+\infty}\biggl[\log\frac
{ f_j(\psi,\theta_j;x)}{f_i(\psi^0,\theta_i^0;x)}\biggr]f_i(\psi^0,\theta
_i^0,x)\,dx
\nonumber
\\[-8pt]
\\[-8pt]
\nonumber
&&\quad\mbox{for } i,j=1,2,\ldots,k+1.
\end{eqnarray}
Note that $J$ may be written $J = J_1 + J_2$, where
%
\begin{equation}\label{eqn2.3}
J_1=\sum_{j=1}^{k+1}\sum_{i=1}^{k+1}\frac{n_{ji}}{n}v(\psi,\theta_j;\psi
^0,\theta_i^0)
\end{equation}
and
%
\begin{eqnarray}\label{eqn2.4}
J_2&=&\frac{1}{n}\sum_{j=1}^{k+1}\sum_{i=n_{j-1}+1}^{n_j}\{\log f_j(\psi
,\theta_j;x_i)-E[\log f_j(\psi,\theta_j;X_i)]\}
\nonumber
\\[-8pt]
\\[-8pt]
\nonumber
&&{}-\frac{1}{n}\sum_{j=1}^{k+1}\sum
_{i=n_{j-1}^0+1}^{n_j^0}\{\log f_j(\psi^0,\theta_j^0;x_i)-E[\log
f_j(\psi^0,\theta_j^0;X_i)]\}.
\end{eqnarray}
Alternatively, we may write
%
\begin{eqnarray}\label{eqn2.5}
J_2&=&\frac{1}{n}\sum_{j=1}^{k+1}\sum_{i=1}^{k+1}\biggl\{\sum_{t\in\tilde
{n}_{ji}}[\log f_j(\psi,\theta_j;x_t)-E(\log f_j(\psi,\theta
_j;X_t))]
\nonumber
\\[-8pt]
\\[-8pt]
\nonumber
&&{}\hspace*{44pt}-\sum_{t\in\tilde{n}_{ji}}[\log f_i(\psi^0,\theta
_i^0;x_t)-E(\log f_i(\psi^0,\theta_i^0,X_t))]\biggr\},
\end{eqnarray}
where $\tilde{n}_{ji}=[n_{j-1}+1,n_j]\cap[n_{i-1}^0+1,n_i^0]$.

Note that $J_1$ is a weighted sum of the negative Kullback--Leibler
distances; it will be shown that $J_2$ approaches $0$ as $n\rightarrow
\infty$.
Also, $v(\psi,\theta_j;\psi^0,\theta_i^0)$ $\leq0$ with equality if
and only if $f_j(\psi,\theta_j;x)=f_i(\psi^0,\theta_i^0;x)$ almost
everywhere (Kullback and Leibler (\citeyear{Kullback51})).

Lemma \ref{lem1} gives a bound for $J_1$.

\begin{lema}\label{lem1}
Under Assumption \textup{\ref{a20}} and Part 1 of Assumption \textup{\ref{a21}}, there exist two
positive constants $C_1>0$ and $C_2>0$ such that, for any $\lambda$ and
$\phi$,
\[
J_1\leq-\max\{C_1\|\lambda-\lambda^0\|_{\infty},C_2\rho(\phi, \phi^0)\},
\]
where $\|\lambda-\lambda^0\|_{\infty}=\max_j|\lambda_j-\lambda_j^0|$ and
$\rho(\phi, \phi^0) =\max_{j}|v(\psi,\theta_j;\psi^0,\theta_j^0)|$.
\end{lema}

Lemma \ref{lem2} describes between-segment properties and within-segment
properties of this model.

\begin{lema}\label{lem2}
Under Part 1 of Assumption \ref{a21}, the following two results follow from
Parts 1 and~2 of Assumption \ref{a22} respectively:
\begin{longlist}[(II)]
\item[(I)] for any $j=1,2,\ldots,k+1$, any $0\leq m_1<m_2\leq n$ and any
positive number $\varepsilon>0$, there exist a constant $A_j$, independent
of $\varepsilon$, and a constant $r<2$, such that
%
\begin{eqnarray}\label{eq6}
&&P_r\Biggl(\max_{m_1\leq s<t\leq m_2, \theta_j\in\Theta_j,\psi\in\Psi}\Biggl|\sum
_{i=s+1}^{t}\{\log f_j(\psi,\theta_j;X_i)-E[\log f_j(\psi,\theta_j;X_i)]\}\Biggr|>\varepsilon\Biggr)
\nonumber
\\[-8pt]
\\[-8pt]
\nonumber
&&\quad \leq A_j\frac
{(m_2-m_1)^r}{\varepsilon^2}.
\end{eqnarray}

\item[(II)] for any $j=1,2,\ldots,k+1$ and any positive number $\varepsilon>0$,
there exist a constant $B_j$, independent of $\varepsilon$, and a constant
$r<2$, such that
%
\begin{eqnarray}\label{eq7}
&&P_r\Biggl(\max_{n_{j-1}^0\leq s<t\leq n_j^0,\psi\in\Psi,\theta_j\in\Theta
_j}\sum_{i=s+1}^t\{[\log f_j(\psi,\theta_j;X_i)
-\log f_j(\psi^0,\theta^0_j;X_i)]
\nonumber
\\[-8pt]
\\[-8pt]
\nonumber
&&{}\hspace*{134pt}-v(\psi,\theta_j;\psi^0,\theta^0_j)\}>\varepsilon\Biggr)\leq B_j\frac
{(n^0_j-n^0_{j-1})^r}{\varepsilon^2}.
\end{eqnarray}
\end{longlist}
\end{lema}

In practical applications, it is useful to have an estimator of $\bar
{i}(\psi^0,\theta^0)$.
Let
\begin{eqnarray*}
\hat{\imath}(\hat{\psi},\hat{\theta})
&=&\pmatrix{
\hat{E}[-\hat{\ell}_{\psi\psi}(\hat{\psi},\hat{\theta});\hat{\phi}]&
\hat{E}[-\hat{\ell}_{\psi\theta}(\hat{\psi},\hat{\theta});\hat{\phi}]\vspace*{3pt}\cr
\hat{E}[-\hat{\ell}_{\psi\theta}(\hat{\psi},\hat{\theta});\hat{\phi
}]^T&\hat{E}[-\hat{\ell}_{\theta\theta}(\hat{\psi},\hat{\theta});\hat
{\phi}]},\\
\hat{E}[-\hat{\ell}_{\psi\psi}(\hat{\psi},\hat{\theta});\hat{\phi
}]&=&\sum_{j=1}^{k+1}\sum_{i=\hat{n}_{j-1}+1}^{\hat{n}_j}\frac
{1}{f_j^2(\hat{\psi},\hat{\theta}_j;x_i)}{f_j}_{\psi}(\hat{\psi},\hat
{\theta}_j;x_i){f_j}_{\psi}^T(\hat{\psi},\hat{\theta}_j;x_i),\\
\hat{E}[-\hat{\ell}_{\psi\theta_j}(\hat{\psi},\hat{\theta});\hat{\phi
}]&=&\sum_{i=\hat{n}_{j-1}+1}^{\hat{n}_j}\frac{1}{f_j^2(\hat{\psi},\hat
{\theta}_j;x_i)}{f_j}_{\psi}(\hat{\psi},\hat{\theta}_j;x_i){f_j}_{\theta
_j}^T(\hat{\psi},\hat{\theta}_j;x_i),\\
\hat{E}[-\hat{\ell}_{\theta_j\theta_j}(\hat{\psi},\hat{\theta});\hat
{\phi}]&=&\sum_{i=\hat{n}_{j-1}+1}^{\hat{n}_j}\frac{1}{f_j^2(\hat{\psi
},\hat{\theta}_j;x_i)}{f_j}_{\theta_j}(\hat{\psi},\hat{\theta
}_j;x_i){f_j}_{\theta_j}^T(\hat{\psi},\hat{\theta}_j;x_i)
\end{eqnarray*}
for $j=1,2,\ldots,k+1$.
Then $\hat{\imath}(\hat{\psi},\hat{\theta})/n$ is a consistent
estimator of $\bar{i}(\psi^0,\theta^0)$.

\section{An example}\label{sec3}

Consider the problem of analyzing the mineral content of a core sample,
which is extensively studied in Chen and Gupta (\citeyear{Chen00}), Chernoff (\citeyear{Chernoff73})
and Srivastava and Worsley (\citeyear{Srivastava86}). In particular, we consider the data
in Chernoff (\citeyear{Chernoff73}) on the mineral content of $12$ minerals in a core
sample measured at $N=53$ equally spaced points. Since some of the
minerals have a very low assay, we follow Chen and Gupta (\citeyear{Chen00}) and
Srivastava and Worsley (\citeyear{Srivastava86}) in analyzing only the $p=5$ variables
$Z_1,Z_8,Z_9,Z_{10}$ and $Z_{12}$ with the highest assays.
Thus, we assume that $(Z_1,Z_8,Z_9,Z_{10},Z_{12})$ has a $5$-variate
normal distribution with a within-segment mean parameter vector and a
variance-covariance matrix that is common to all segments.
The analyses of Chen and Gupta (\citeyear{Chen00}), Chernoff (\citeyear{Chernoff73}) and Srivastava
and Worsley (\citeyear{Srivastava86}) suggest that there are $5$ change points of the mean
vector and, hence, we make
that assumption here.

The estimates of $5$ change points, within-segment parameters of mean
vectors and common parameter of variance-covariance matrix were
computed using maximum likelihood. The estimated change points are
$7,20,24,32$ and $41$, which are different from those estimated change
points by Chen and Gupta (\citeyear{Chen00}), Chernoff (\citeyear{Chernoff73}) and Srivastava and
Worsley~(\citeyear{Srivastava86}), and are more reasonable. This is because Chen and Gupta
(\citeyear{Chen00}), Chernoff (\citeyear{Chernoff73}) and Srivastava and Worsley~(\citeyear{Srivastava86}) use the
binary segmentation procedures which detect multiple change points one
by one, not simultaneously, whereas the method in this paper
simultaneously estimates multiple change points.
The estimated six within-segment mean vectors are in the
following.  They are arranged according to the order of from left to right.
For example, the two vectors on the first line are, respectively, the first
and second within-segment mean vectors.
\begin{eqnarray*}
&(287.14,58.57,25.71,240.00,422.86),\qquad
(277.31,144.61,24.69,306.15,274.62),&\\
&(321.25,502.50,150.00,620.00,217.50),\qquad(397.50,635.00,428.75,625.00,4.38),&
\\
&(470.00,188.89,214.44,255.56,108.89),\qquad
(425.0,155.92,183.42,320.0,333.33).&
\end{eqnarray*}
The estimated common variance-covariance matrix is
\[
\pmatrix{
1485.71&-966.03&569.41&-421.41&-590.87\cr
-966.03&8523.65&4649.95&5982.95&1054.22\cr
569.41&4649.95&8767.11&4434.76&736.33\cr
-421.41&5982.95&4434.76&8768.49&780.03\cr
-590.87&1054.22&736.33&780.03&3193.37}.
\]

\section{Discussion}\label{sec4}
This paper establishes the consistency of maximum
likelihood estimators of the parameters of a general class of multiple
change-point models and gives the asymptotic distribution of the
parameters of the within-segment distributions.
The required regularity conditions are relatively weak and are
generally satisfied by exponential family distributions.

Some important problems in the analysis of multiple change-point models
were not considered here. One is that the asymptotic distribution of the
maximum likelihood estimator of the vector of change points was not
considered. The reason for this is that
the methods used to determine this asymptotic distribution are quite
different from the methods used to
establish the consistency of the maximum likelihood estimator; see, for
example, Hinkley (\citeyear{Hinkley70}) for a treatment of this
problem in a single change-point model. Thus, this is essentially a
separate research topic. However, the asymptotic properties obtained in
this paper are necessary for the establishment of the asymptotic
distribution of the maximum likelihood estimator of the vector of
change points in this model. This will be a subject of future work.

Another important problem is to extend the results of this paper to the
case in which the number of change points is not known and must be
determined from the data. Clearly, a likelihood-based approach to this
problem will require
an understanding of the properties of maximum likelihood estimators in
the model in which the number of change points is known. Thus, the
results of the present paper can be considered as a first step toward
the development of a
likelihood-based methodology that can be used to determine
simultaneously the number and location of the change points. This is
also a topic of future research.
\appendix
\section*{\texorpdfstring{Appendix A: The consistency assumption of Hinkley (\citeyear{Hinkley72})}
{Appendix A: The consistency assumption of Hinkley (1972)}}\label{appA}

Consider a change-point model with a single change point, $n_1^0$, and
suppose that there are no common parameters in the model.
In Hinkley (\citeyear{Hinkley72}), it is shown that ${\hat n}_1$, the maximum
likelihood estimator
of $n_1^0$, satisfies ${\hat n}_1 = n_1^0 + \mathrm{O}_p(1)$ under the
condition
\setcounter{equation}{0}
\begin{equation}\label{eqa1}
\sup_{\theta_1}\sum_{i=n_1^0 +1}^{n_1^0 + m} \{ \log f_1(X_i; \theta_1)
- \log f_2(X_i; \theta_2^0) \} \to-\infty
\end{equation}
with probability $1$ as $m\to\infty$, which was described as a
``consistency assumption''. Note that the random variables in the sum
$X_{n_1^0 + 1}, \ldots, X_{n_1^0 + m}$ are drawn from the distribution
with density $f_2$.

Suppose that
\[
{1\over m} \sum_{i=n_1^0 +1}^{n_1^0 + m} \{ \log f_1(X_i; \theta_1) -
\log f_2(X_i; \theta_2^0) \}
\]
converges to
\[
E\biggl\{ \log{f_1(X; \theta_1) \over f_2(X; \theta_2^0)}\biggr\}
\]
as $m \to\infty$, uniformly in $\theta_1$, where $X$ is distributed
according to the distribution with density $f_2(\cdot; \theta_2^0)$.
Equation (\ref{eqa1}) then holds, provided that
\[
\sup_{\theta_1}E\biggl\{ \log{f_1(X; \theta_1) \over f_2(X; \theta_2^0)}\biggr\} <
0;
\]
note that, by properties of the Kullback--Leibler distance and
Assumption \ref{a20},
\[
E\biggl\{ \log{f_1(X; \theta_1) \over f_2(X; \theta_2^0)}\biggr\}<0
\]
for each $\theta_1$.

Thus, condition (\ref{eqa1}) fails whenever the distribution corresponding
to the density $f_2(\cdot; \theta_2^0)$ is in the closure of the
set of distributions corresponding to densities of the form
$f_1(\cdot; \theta_1)$, in a certain sense.

One such case occurs if $f_1$ and $f_2$ have the same parametric form
with parameters $\theta_1^0, \theta_2^0$, respectively, satisfying
$\theta_1^0 \not= \theta_2^0$. For instance, suppose that the random
variables in the first segment are normally distributed with mean
$\theta_1^0$ and standard deviation $1$ and the random variables in the
second segment are normally distributed with mean $\theta_2^0$ and
standard deviation $1$. Then
\[
\sup_{\theta_1}\sum_{i=n_1^0 +1}^{n_1^0 + m} \{ \log f_1(X_i; \theta_1)
- \log f_2(X_i; \theta_2^0) \} = {m \over2} ({\bar X}_m - \theta
_2^0)^2 ,
\]
where
\[
{\bar X}_m = \frac{1}{m}\sum_{i = n_1^0 + 1}^{n_1^0+m} X_i
\]
is normally distributed with mean $\theta_2^0$ and variance $1/m$.
Clearly, (\ref{eqa1}) does not hold in this case.

A similar situation occurs when the distribution with density $f_2(\cdot
; \theta_2^0)$ can be viewed as a~limit of the distributions
with densities $f_1(\cdot; \theta_1)$. For instance, suppose that
$f_1$ is the density of a~Weibull distribution with rate parameter
$\beta$ and shape parameter $\alpha$, $\theta_1 = (\alpha, \beta)$,
$\beta\not=1$, and $f_2$ is the density of an exponential distribution
with rate parameter $\theta_2$.

In this appendix, we show that this is a strong assumption that is not
generally satisfied by otherwise well-behaved
models. For instance, suppose that $f_1$ and $f_2$ have the same
functional form and that the
difference between the two distributions is due to the fact that $\theta
_1^0 \not= \theta_2^0$. Again, (\ref{eqa1}) will not hold.

Thus, the consistency condition used in Hinkley (\citeyear{Hinkley72}) is too strong
for the general model considered here.

\section*{Appendix B: Technical details}\label{appB}

\begin{pf*}{Proof of Lemma \protect\ref{lem1}}
We first need to prepare some results which are to be used in this
proof. For $i=1,2,\ldots,k$, let us define
\[
g_i(\alpha,\phi^0)=\sup_{1\leq j\leq k+1}\sup_{\theta_j\in\Theta_j}\sup
_{\psi\in\Psi}[\alpha v(\psi,\theta_j;\psi^0,\theta_{i+1}^0)+(1-\alpha
)v(\psi,\theta_j;\psi^0,\theta_{i}^0)],
\]
where $0\leq\alpha\leq1$.
We then have that $g_i(0,\phi^0)=g_i(1,\phi^0)=0$ for $i=1,2,\ldots,k$.
It is straightforward to show that $g_i(\alpha,\phi^0)$ is a convex
function with respect to $\alpha$ for any $i=1,2,\ldots,k$.

Let $G_i(\phi^0)=2g_i(1/2,\phi^0)$. Because $\alpha=2\alpha
(1/2)+(1-2\alpha)0$ for $0\leq\alpha\leq1/2$, convexity of $g_i(\alpha
,\phi^0)$ gives that
\[
g_i(\alpha,\phi^0)\leq2\alpha g_i(1/2,\phi^0)=\alpha G_i(\phi^0)\qquad\mbox
{for }i=1,2,\ldots,k.
\]
Noting that
\[
g_i(1/2,\phi^0)=\frac{1}{2}\sup_{1\leq j\leq k+1}\sup_{\theta_j\in\Theta
_j}\sup_{\psi\in\Psi}[ v(\psi,\theta_j;\psi^0,\theta_{i+1}^0)+v(\psi
,\theta_j;\psi^0,\theta_{i}^0)],
\]
it follows from Assumption \ref{a20} that $G_i(\phi^0)<0$. If we let $\bar
{G}(\phi^0)=\max_{1\leq i\leq k}G_i(\phi^0)$, then $\bar{G}(\phi^0)<0$.

Let $\Delta_{\lambda}^0=\min_{1\leq j\leq k-1}|\lambda_{j+1}^0-\lambda
_j^0|$. Consider a change-point fraction configuration $\lambda$ such
that $\|\lambda-\lambda^0\|_{\infty}\leq\Delta_{\lambda}^0/4$. For any
$j=1,2,\ldots,k$, there are two cases: a candidate change-point fraction
$\lambda_j$ may be on the left or the right of the true change-point
fraction $\lambda_j^0$.

For any $j$ with $\lambda_j$ on the right of $\lambda_j^0$, we have
that $\lambda_{j-1}\leq\lambda_j^0\leq\lambda_j$. Then
\[
J_1\leq\frac{n_{j,j+1}}{n}v(\psi,\theta_j;\psi^0,\theta_{j+1}^0)+\frac
{n_{jj}}{n}v(\psi,\theta_j;\psi^0,\theta_j^0).
\]
If we define $\alpha_{j,j+1}=n_{j,j+1}/(n_{j,j+1}+n_{jj})$, then the
case $\|\lambda-\lambda^0\|_{\infty}\leq\Delta_{\lambda}^0/4$ gives
that $\alpha_{j,j+1}\leq\frac{1}{2}$ and
\begin{eqnarray*}
J_1&\leq&\frac{n_{j,j+1}+n_{jj}}{n}[\alpha_{j,j+1}v(\psi
,\theta_j;\psi^0,\theta_{j+1}^0)+(1-\alpha_{j,j+1})v(\psi,\theta_j;\psi
^0,\theta_j^0)]\\
&\leq&\frac{n_{j,j+1}}{n}G_j(\phi^0)\leq(\lambda_j-\lambda_j^0)\bar
{G}(\phi^0).
\end{eqnarray*}

For any $j$ with $\lambda_j$ on the left of $\lambda_j^0$, we have that
$\lambda_j\leq\lambda_j^0\leq\lambda_{j+1}$.\vspace*{-1pt} Similarly, we define
$\alpha_{j,j-1}=n_{j,j-1}/(n_{j,j-1}+n_{jj})$. Using the fact that
$\alpha_{j,j-1}\leq\frac{1}{2}$, it similarly gives that $J_1\leq
(\lambda_j^0-\lambda_j)\bar{G}(\phi^0)$.

Therefore, if $\|\lambda-\lambda^0\|_{\infty}\leq\Delta_{\lambda}^0/4$,
then we obtain that $J_1\leq\|\lambda-\lambda^0\|_{\infty}\bar{G}(\phi
^0)$. On the other hand,
\[
J_1\leq\min_{1\leq j\leq k+1}v(\psi,\theta_j;\psi^0,\theta_j^0)\frac
{n_{jj}}{n}=-\max_{1\leq j\leq k+1}|v(\psi,\theta_j;\psi^0,\theta
_j^0)|\frac{n_{jj}}{n}.
\]
We have $n_{jj}/n\geq\Delta_{\lambda}^0/2$ for any $j$, so
\[
J_1\leq-\frac{1}{2}\Delta_{\lambda}^0\sup_{1\leq j\leq k+1}|v(\psi
,\theta_j;\psi^0,\theta_j^0)|=-\frac{1}{2}\Delta_{\lambda}^0\rho(\phi
,\phi^0).
\]

Now, consider the other case of a change-point fraction configuration
$\lambda$, where $\|\lambda-\lambda^0\|_{\infty}>\Delta_{\lambda}^0/4$.
It is clear that there exists a pair of integers $(i,j)$ such that
$n_{ij}\geq n\Delta_{\lambda}^0/4$, $n_{i,j+1}\geq n\Delta_{\lambda
}^0/4$ and $n_{ij}\geq n_{i,j+1}$. Let $\alpha
_{i,j+1}=n_{i,j+1}/(n_{i,j+1}+n_{ij})$. For any $\phi$, we have that
\begin{eqnarray*}
J_1&\leq&\frac{n_{i,j+1}+n_{ij}}{n}[\alpha_{i,j+1}v(\psi,\theta_i;\psi
^0,\theta_{j+1}^0)+(1-\alpha_{i,j+1})v(\psi,\theta_i;\psi^0,\theta
_j^0)]\\
&\leq&\frac{n_{i,j+1}+n_{ij}}{n}\min(\alpha_{i,j+1},1-\alpha
_{i,j+1})\bar{G}(\phi^0)\\
&\leq&\frac{\Delta_{\lambda}^0}{2}\min\biggl(\frac
{n_{i,j+1}}{n},\frac{n_{ij}}{n}\biggr)\bar{G}(\phi^0)\\
&\leq&\frac{1}{2}\biggl(\frac{{\Delta_{\lambda}^0}}{2}\biggr)^2\bar
{G}(\phi^0).
\end{eqnarray*}

Combining the results from the two cases of $\|\lambda-\lambda^0\|
_{\infty}\leq\Delta_{\lambda}^0/4$ and $\|\lambda-\lambda^0\|_{\infty
}>\Delta_{\lambda}^0/4$, it follows that
\[
J_1\leq\bar{G}(\phi^0)\min\biggl(\frac{1}{2}\biggl(\frac{\Delta_{\lambda
}^0}{2}\biggr)^2,\|\lambda-\lambda^0\|_{\infty}\biggr)\leq\frac{1}{2}\biggl(\frac{\Delta
_{\lambda}^0}{2}\biggr)^2\bar{G}(\phi^0)\|\lambda-\lambda^0\|_{\infty}
\]
and

\renewcommand{\theequation}{B.\arabic{equation}}
\setcounter{equation}{0}
\begin{equation}\label{eqb1}
J_1\leq\frac{\Delta_{\lambda}^0}{2}\max\biggl[-\rho(\phi,\phi^0),
\frac{\Delta_{\lambda}^0}{4}\bar{G}(\phi^0)\biggr]
\leq-\frac{\Delta_{\lambda}^0}{2}\min\biggl[\rho(\phi,\phi^0),-\frac{\Delta
_{\lambda}^0}{4}\bar{G}(\phi^0)\\\biggr] .
\end{equation}

Note that (\ref{eqb1}) can be simplified. If we define
\[
\varrho(
\phi,\phi^0)=\max_{1\leq j\leq k+1}\sup_{\theta_j\in\Theta_j}\sup_{\psi
\in\Psi}|v(\psi,\theta_j;\psi^0,\theta_j^0)|,
\]
then we have that $\rho(\phi,\phi^0)/\varrho(\phi,\phi^0)\leq1$. It
follows from inequality (\ref{eqb1}) that
\[
J_1\leq-\frac{\Delta_{\lambda}^0}{2}\varrho(\phi,\phi^0)\min\biggl[\frac{\rho
(\phi,\phi^0)}{\varrho(\phi,\phi^0)},-\frac{\Delta_{\lambda}^0}{4}\bar
{G}(\phi^0)/\varrho(\phi,\phi^0)\biggr].
\]
If $-(\Delta_{\lambda}^0/4)\bar{G}(\phi^0)/\varrho(\phi,\phi^0)\leq1$,
then we have that
\[
J_1\leq(\Delta_{\lambda}^0/2)^2\bigl(\rho(\phi,\phi^0)/\varrho(\phi,\phi
^0)\bigr)\bigl(\bar{G}(\phi^0)/2\bigr).
\]
If $-(\Delta_{\lambda}^0/4)\bar{G}(\phi^0)/\varrho(\phi,\phi^0)>1$,
then $J_1\leq-(\Delta_{\lambda}^0/2)\rho(\phi,\phi^0)$. Letting
\[
C_2=\min\{(\Delta_{\lambda}^0/2)^2|\bar{G}(\phi^0)|/(2\varrho(\phi,\phi
^0)),\Delta_{\lambda}^0/2\},
\]
inequality (\ref{eqb1}) gives that $J_1\leq-C_2\rho(\phi,\phi^0)$.

Setting $C_1=(\Delta_{\lambda}^0/2)^2|\bar{G}(\phi^0)|/2$, we finally
have that
\[
J_1\leq-\max\{C_1\|\lambda-\lambda^0\|_{\infty},C_2\rho(\phi,\phi^0)\},
\]
which concludes the proof.
\end{pf*}

\begin{pf*}{Proof of Lemma \protect\ref{lem2}}
With Part 1 of Assumption \ref{a22} in mind, equation (\ref{eq6}) can be achieved by
induction with respect to $m_2$. The induction method is similar to the
one used in M\'{o}ricz, Serfling and Stout~(\citeyear{Moricz82}), so its
proof is omitted here. Using Part 2 of Assumption \ref{a22}, equation (\ref{eq7}) can
be proven similarly by the same induction method.
\end{pf*}
\begin{pf*}{Proof of Theorem \protect\ref{thm1}}
Let
\begin{eqnarray*}
\Lambda_{\delta}&=&\{\lambda\in\Lambda\dvt\|\lambda-\lambda^0\|_{\infty
}>\delta\},\qquad
 \Phi_{\delta}=\{\phi\in\Phi\dvt \rho(\phi,\phi^0)>\delta\},
\\
\Phi&=&\Theta_1\times\Theta_2\times\cdots\times\Theta_{k+1}\times\Psi,
\\
\Lambda&=&\{(\lambda_1,\lambda_2,\ldots,\lambda_k)|\lambda_j=n_j/n,
j=1,2,\ldots,k;\\
&&{}\hspace*{6pt} 0<n_1<n_2<\cdots<n_k<n\}.
\end{eqnarray*}
Then, for any $\delta>0$, it follows from Lemma \ref{lem1} that
\begin{eqnarray}
-\max_{\lambda\in\Lambda_{\delta},\phi\in\Phi}J_1\geq C_1\delta\quad \mbox{and}\quad
-\max_{\phi\in\Phi_{\delta},\lambda\in\Lambda}J_1\geq
C_2\delta.\nonumber
\end{eqnarray}
Therefore, we obtain that
\begin{eqnarray*}
&&P_r(\|\hat{\lambda}-\lambda^0\|_{\infty}>\delta)\\
&&\quad\leq P_r\Bigl(\max_{\lambda
\in\Lambda_{\delta},\phi\in\Phi}J>0\Bigr) \leq P_r\Bigl(\max_{\lambda\in\Lambda_{\delta},\phi\in\Phi}J_2>-\max
_{\lambda\in\Lambda_{\delta},\phi\in\Phi}J_1\Bigr) \leq P_r\Bigl(\max_{\lambda\in\Lambda_{\delta},\phi\in\Phi}|J_2|>C_1\delta\Bigr)\\
&&\quad \leq P_r\Biggl(\max_{\lambda\in\Lambda_{\delta},\phi\in\Phi}\sum
_{j=1}^{k+1}\frac{1}{n}\Biggl|\sum_{i=n_{j-1}+1}^{n_j}\{\log f_j(\psi,\theta
_j;X_i)-E[\log f_j(\psi,\theta_j;X_i)]\}\Biggr|>\frac{C_1}{2}\delta\Biggr)
\\
&&\qquad{}+P_r\Biggl(\sum_{j=1}^{k+1}\frac{1}{n}\Biggl|\sum
_{i=n_{j-1}^0+1}^{n_j^0}\{\log f_j(\psi^0,\theta_j^0;X_i)-E[\log
f_j(\psi^0,\theta_j^0;X_i)]\Biggr|>\frac{C_1}{2}\delta\Biggr)\\
&&\quad\leq\sum_{j=1}^{k+1}P_r\Biggl(\max_{0\leq n_{j-1}<n_j\leq n,\theta_j\in
\Theta_j,\psi\in\Psi}\frac{1}{n}\Biggl|\sum_{i=n_{j-1}+1}^{n_j}\{\log f_j(\psi
,\theta_j;X_i)-E[\log f_j(\psi,\theta_j;X_i)]\}\Biggr|\\
&&\quad\hspace*{43pt}>\frac{C_1\delta
}{2(k+1)}\Biggr)\\
&&\qquad{}+\sum_{j=1}^{k+1}P_r\Biggl(\frac{1}{n}\Biggl|\sum
_{i=n_{j-1}^0+1}^{n_j^0}\{\log f_j(\psi^0,\theta_j^0;X_i)-E[\log
f_j(\psi^0,\theta_j^0;X_i)]\}\Biggr|>\frac{C_1\delta}{2(k+1)}\Biggr) .
\end{eqnarray*}

It follows from Lemma \ref{lem2} that
\[
P_r(\|\hat{\lambda}-\lambda^0\|_{\infty}>\delta)\leq2\biggl[\frac
{2(k+1)}{C_1\delta}\biggr]^2\Biggl(\sum_{j=1}^{k+1}A_j\Biggr)n^{r-2}\rightarrow0\qquad \mbox
{as } n\rightarrow+\infty,
\]
noting that $r<2$.

For $\hat{\phi}$, we similarly obtain that
\begin{eqnarray*}
&&P_r\bigl(\rho(\hat{\phi},\phi^0)>\delta\bigr)\\
&&\quad\leq P_r\Bigl(\max_{\lambda\in\Lambda
,\phi\in\Phi_{\delta}}J>0\Bigr)\\
&&\quad \leq\sum_{j=1}^{k+1}P_r\Biggl(\max_{0\leq n_{j-1}<n_j\leq n,\theta_j\in
\Theta_j,\psi\in\Psi}\frac{1}{n}\Biggl|\sum_{i=n_{j-1}+1}^{n_j}\{\log f_j(\psi
,\theta_j;X_i)-E[\log f_j(\psi,\theta_j;X_i)]\}\Biggr|\\
&&{}\hspace*{51pt}>\frac{C_2\delta
}{2(k+1)}\Biggr)\\
&&\qquad{}+\sum_{j=1}^{k+1}P_r\Biggl(\frac{1}{n}\Biggl|\sum
_{i=n_{j-1}^0+1}^{n_j^0}\{\log f_j(\psi^0,\theta_j^0;X_i)-E[\log
f_j(\psi^0,\theta_j^0;X_i)]\}\Biggr|>\frac{C_2\delta}{2(k+1)}\Biggr).
\end{eqnarray*}
Similarly, Lemma \ref{lem2} shows that $P_r(\rho(\hat{\phi},\phi^0)>\delta
)\rightarrow0\ \mbox{as}\ n\rightarrow+\infty$. Noting the fact that
$v(\psi,\theta_j;\psi^0,\theta^0_j)=0$ if and only if $\psi=\psi^0$ and
$\theta_j=\theta^0_j$,\vspace*{1pt} it follows that $\hat{\psi}\rightarrow_p\psi^0$
and $\hat{\theta}_j\rightarrow_p\theta^0_j$ for $j=1,2,\ldots,k+1$,
which completes the proof.
\end{pf*}

\begin{pf*}{Proof of Theorem \protect\ref{thm2}}
Let us first define
\[
\Lambda_{\delta,n}=\{\lambda\in\Lambda\dvt n\|\lambda-\lambda^0\|_{\infty
}>\delta\}
\]
for any $\delta>0$. Because of the consistency of $\hat{\lambda}$, we
need to consider only those terms whose observations are in $\tilde
{n}_{j,j-1}$, $\tilde{n}_{j,j}$ and $\tilde{n}_{j,j+1}$ for all $j$ in
equation (\ref{eqn2.5}). Therefore, we have
\begin{eqnarray*}\label{formu3}
&&P_r(n\|\hat{\lambda}-\lambda^0\|_{\infty}>\delta)\\
&&\quad \leq\sum_{j=1}^{k+1}P_r\Biggl(\max_{\lambda\in\Lambda_{\delta,n},\phi\in\Phi
}\Biggl\{\frac{1}{n}\sum_{t\in\tilde{n}_{jj}}[\log f_j(\psi,\theta
_j;X_t)-E(\log f_j(\psi,\theta_j;X_t))]\\
&&\qquad{}\hspace*{78pt}-\frac{1}{n}\sum_{t\in\tilde{n}_{jj}}[\log f_j(\psi
^0,\theta^0_j;X_t)-E(\log f_j(\psi^0,\theta^0_j;X_t))]\\
&&\qquad{}\hspace*{78pt}+\frac{1}{3(k+1)}J_1\Biggr\}>0\Biggr)\\
&&\qquad{}+\sum_{j=2}^{k+1}P_r\Biggl(\max_{\lambda\in\Lambda_{\delta
,n},\phi\in\Phi}\Biggl\{\frac{1}{n}\sum_{t\in\tilde{n}_{j,j-1}}[\log f_j(\psi
,\theta_j;X_t)-E(\log f_j(\psi,\theta_j;X_t))]\\
&&\qquad{}\hspace*{92pt}-\frac{1}{n}\sum_{t\in\tilde{n}_{j,j-1}}[\log f_{j-1}(\psi
^0,\theta^0_{j-1};X_t)-E(\log f_{j-1}(\psi^0,\theta^0_{j-1};X_t))]\\
&&\qquad{}\hspace*{92pt}+\frac
{1}{3k}J_1\Biggr\}>0\Biggr)\\
&&\qquad{}+\sum_{j=1}^{k}P_r\Biggl(\max_{\lambda\in\Lambda_{\delta,n},\phi
\in\Phi}\Biggl\{\frac{1}{n}\sum_{t\in\tilde{n}_{j,j+1}}[\log f_j(\psi,\theta
_j;X_t)-E(\log f_j(\psi,\theta_j;X_t))]\\
&&\qquad{}\hspace*{91pt}-\frac{1}{n}\sum_{t\in\tilde{n}_{j,j+1}}[\log f_{j+1}(\psi
^0,\theta^0_{j+1};X_t)-E(\log f_{j+1}(\psi^0,\theta^0_{j+1};X_t))]\\
&&\qquad{}\hspace*{92pt}+\frac
{1}{3k}J_1\Biggr\}>0\Biggr)\\
&&\quad \equiv\sum_{j=1}^{k+1}I_{1j}+\sum_{j=2}^{k+1}I_{2j}+\sum
_{j=1}^{k}I_{3j}.
\end{eqnarray*}

First, consider the probability formulas $I_{1j}$ in the above equation
for any $j=1,2,\ldots,k+1$. The consistency of $\hat{\lambda}$ allows us
to restrict our attention to the case $n_{jj}>\frac
{1}{2}(n^0_j-n^0_{j-1})$. For this case, we have that
\[
J_1\leq\frac{n^0_j-n^0_{j-1}}{2n}v(\psi,\theta_j;
\psi^0,\theta_j^0) .
\]

Therefore, we obtain that
\begin{eqnarray*}
I_{1j}&\leq& P_r\biggl(\sum_{t\in\tilde{n}_{jj}^*}\{[\log f_j(\psi^*,\theta
^*_j;x_t)-\log f_j(\psi^0,\theta^0_j;x_t)]\\
&&\qquad{}\hspace*{16pt}-v(\psi^*,\theta^*_j;\psi
^0,\theta^0_j)\}>\frac{n^0_j-n^0_{j-1}}{6(k+1)} |v(\psi^*,\theta
_j^*;\psi^0,\theta^0_j)|\biggr)\\
&\leq& P_r\Biggl(\max_{n_{j-1}^0\leq s<t\leq n_j^0,\psi\in\Psi
,\theta_j\in\Theta_j}\sum_{i=s+1}^t\{[\log f_j(\psi,\theta
_j;X_t)-\log f_j(\psi^0,\theta^0_j;X_t)]\\
&&\qquad{}\hspace*{112pt}-v(\psi,\theta_j;\psi^0,\theta^0_j)\}>\frac
{E}{6(k+1)}(n^0_j-n^0_{j-1})\Biggr),
\end{eqnarray*}
where $\tilde{n}^*_{jj}$, $\psi^*$, $\theta_j^*$ and $\lambda^*$ are,
respectively, the maximizing values of\vspace*{1pt} $\tilde{n}_{jj}$, $\psi$, $\theta
_j$ and $\lambda$
obtained through the maximization. Equation (\ref{eq7}) of Lemma \ref{lem2} can then
be applied to show that $I_{1j}\rightarrow0$ as $n,\delta\rightarrow
\infty$.

Next, consider the probability formula $I_{2j}$ for any $j=2,\ldots
,k+1$. In this case, $\lambda_{j-1}<\lambda^0_{j-1}$. We have that
\begin{eqnarray*}
I_{2j} &\leq& P_r\biggl(\max_{\lambda\in\Lambda_{\delta,n},\phi\in\Phi}\biggl\{
\frac{1}{n}\sum_{t\in\tilde{n}_{j,j-1}}[\log f_j(\psi,\theta
_j;X_t)-E(\log f_j(\psi,\theta_j;X_t))]+\frac{1}{6k}J_1\biggr\}>0\biggr)\\
&&{}\hspace*{18pt}+P_r\biggl(\max_{\lambda\in\Lambda_{\delta,n},\phi\in\Phi}\biggl\{
-\frac{1}{n}\sum_{t\in\tilde{n}_{j,j-1}}[\log f_{j-1}(\psi,\theta
_{j-1};X_t)-E(\log f_{j-1}(\psi,\theta_{j-1};X_t))]\\
&&\qquad{}\hspace*{74pt}+\frac{1}{6k}J_1\biggr\}>0\biggr)\\
&\equiv& I_{2j}^{(1)}+I_{2j}^{(2)}.
\end{eqnarray*}

$I_{2j}^{(1)}$ and $I_{2j}^{(2)}$\vspace*{1pt} can be handled in the same way, so we
just show how to handle $I_{2j}^{(1)}$. Only two cases have to be considered.

If $n^0_{j-1}-n_{j-1}\leq\delta$, then
\begin{eqnarray*}
&&I_{2j}^{(1)}\leq P_r\Biggl(\max_{n_{j-1}\leq s<t\leq n^0_{j-1},\theta_j\in
\Theta_j,\psi\in\Psi}\Biggl|\sum_{i=s+1}^{t}[\log f_j(\psi,\theta_j;X_i)
-E(\log f_j(\psi,\theta_j;X_i))]\Biggr|\\
&&{}\quad\hspace*{32pt}>\frac{C_1\delta}{6k}\Biggr).
\end{eqnarray*}
Equation (\ref{eq6}) of Lemma \ref{lem2} gives that $I_{2j}^{(1)}\rightarrow0$ as
$n,\delta\rightarrow+\infty$.

If $n^0_{j-1}-n_{j-1}>\delta$ for the other case, then $J_1\leq
-C_1(n^0_{j-1}-n_{j-1})/n$. Therefore, we obtain that
\begin{eqnarray*}
I^{(1)}_{2j}&\leq& P_r\Biggl(\max_{n_{j-1}\leq s<t\leq n^0_{j-1},\theta_j\in
\Theta_j,\psi\in\Psi}\frac{n^0_{j-1}-n_{j-1}}{n}\\
&&{}\hspace*{106pt}\quad \times\Biggl(\frac
{1}{n^0_{j-1}-n_{j-1}}\sum_{i=s+1}^{t}[\log f_j(\psi,\theta_j;X_i)\\
&&\qquad{}\hspace*{193pt}-E(\log f_j(\psi,\theta_j;X_i))]\\
&&\qquad{}\hspace*{113pt}-\frac{C_1}{6k}\Biggr)>0\Biggr)\\
&\leq& P_r\Biggl(\max_{n_{j-1}\leq s<t\leq n^0_{j-1},\theta
_j\in\Theta_j,\psi\in\Psi}\Biggl|\sum_{i=s+1}^{t}[\log f_j(\psi,\theta
_j;X_i)-E(\log f_j(\psi,\theta_j;X_i))]\Biggr|\\
&&\quad{}\hspace*{2pt}>\frac
{C_1}{6k}(n^0_{j-1}-n_{j-1})\Biggr),
\end{eqnarray*}
which converges to zero as $n,\delta\rightarrow0$, by equation (\ref{eq6}) of
Lemma \ref{lem2}.

$I_{3j}$ can be handled in the same way as $I_{2j}$. Therefore, Theorem
\ref{thm2} is proved.
\end{pf*}
\begin{pf*}{Proof of Theorem \protect\ref{thm3}}
We first have the expansion
\[
\hat{\ell}_{\phi}(\hat{\psi},\hat{\theta})-\hat{\ell}_{\phi}(\psi
^0,\theta^0)=[\hat{\ell}_{\phi\phi}(\psi^0,\theta^0)+\mathrm{o}_p(n)](\hat{\phi
}-\phi^0).
\]
The fact that $\hat{\ell}_{\phi}(\hat{\psi},\hat{\theta})=0$ then gives that
\[
\sqrt{n}(\hat{\phi}-\phi^0)=\biggl[-\frac{1}{n}\hat{\ell}_{\phi\phi}(\psi
^0,\theta^0)+\mathrm{o}_p(1)\biggr]^{-1}\frac{\hat{\ell}_{\phi}(\psi^0,\theta^0)}{\sqrt{n}}.
\]

Now, consider the limit of $\hat{\ell}_{\phi}(\psi^0,\theta^0)/\sqrt
{n}$. We have that
\begin{eqnarray}
\frac{1}{\sqrt{n}}\hat{\ell}_{\phi}(\psi^0,\theta^0)
=\frac{1}{\sqrt{n}}[\hat{\ell}_{\phi}(\psi^0,\theta^0)-\ell^0_{\phi
}(\psi^0,\theta^0)]
+\frac{1}{\sqrt{n}}\ell^0_{\phi}(\psi^0,\theta^0).\nonumber
\end{eqnarray}

Because of the consistency of $\hat{\lambda}$, we can assume that
$n^0_{j-1}<\hat{n}_j<n^0_{j+1}$ for $j=1,2,\ldots,k$. It is then
straightforward to obtain that
\begin{eqnarray*}
&&\frac{1}{\sqrt{n}}[\hat{\ell}_{\phi}(\psi^0,\theta^0)-\ell^0_{\phi
}(\psi^0,\theta^0)]\\
&&\quad=\frac{1}{\sqrt{n}}\sum_{j=1}^{k+1}\bigl[\hat{\ell}_{\phi
}^{(j)}(\psi^0,\theta_j^0)-\ell_{\phi}^{(j)}(\psi^0,\theta
_j^0)\bigr]\nonumber\\
&&\quad =\frac{1}{\sqrt{n}}\sum_{j=1}^{k+1}\Biggl\{
I(\hat{n}_j\geq n^0_j,\hat{n}_{j-1}\geq n^0_{j-1})\\
&&{}\hspace*{56pt}\times\Biggl[\sum_{i=n^0_j+1}^{\hat{n}_j}\frac{\partial}{\partial\phi
}\log f_j(\psi^0,\theta^0_j;X_i)
-\sum_{i=n^0_{j-1}+1}^{\hat{n}_{j-1}}\frac{\partial}{\partial\phi}\log
f_j(\psi^0,\theta^0_j;X_i)\Biggr]\\
&&{}\hspace*{58pt}+I(\hat{n}_j\geq n^0_j,\hat{n}_{j-1}< n^0_{j-1})\\
&&{}\hspace*{67pt}\times\Biggl[\sum_{i=n^0_j+1}^{\hat{n}_j}\frac{\partial}{\partial\phi
}\log f_j(\psi^0,\theta^0_j;X_i)+\sum_{i=\hat
{n}_{j-1}+1}^{n^0_{j-1}}\frac{\partial}{\partial\phi}\log f_j(\psi
^0,\theta^0_j;X_i)\Biggr]\\
&&{}\hspace*{58pt}+I(\hat{n}_j<n^0_j,\hat{n}_{j-1}\geq n^0_{j-1})\\
&&{}\hspace*{66pt}\times\Biggl[-\sum_{i=\hat{n}_j+1}^{n^0_j}\frac{\partial}{\partial\phi
}\log f_j(\psi^0,\theta^0_j;X_i)-\sum_{i=n^0_{j-1}+1}^{\hat
{n}_{j-1}}\frac{\partial}{\partial\phi}\log f_j(\psi^0,\theta
^0_j;X_i)\Biggr]\\
&&{}\hspace*{58pt}+I(\hat{n}_j<n^0_j,\hat{n}_{j-1}< n^0_{j-1})\\
&&{}\hspace*{68pt}\times\Biggl[-\sum_{i=\hat{n}_j+1}^{n^0_j}\frac{\partial}{\partial\phi
}\log f_j(\psi^0,\theta^0_j;X_i)+\sum_{i=\hat
{n}_{j-1}+1}^{n^0_{j-1}}\frac{\partial}{\partial\phi}\log f_j(\psi
^0,\theta^0_j;X_i)\Biggr]\Biggr\}.
\end{eqnarray*}

It follows from Theorem \ref{thm2} that
\[
\frac{1}{\sqrt{n}}[\hat{\ell}_{\phi}(\psi^0,\theta^0)-\ell^0_{\phi}(\psi
^0,\theta^0)]=\frac{1}{\sqrt{n}}\mathrm{O}_p(1),
\]
which converges to zero in probability as $n\to\infty$.

Since
\[
\frac{1}{\sqrt{n}}\ell^0_{\phi}(\psi^0,\theta^0)\stackrel{\mathcal
{D}}\rightarrow N_{d+d_1+d_2+\cdots+d_{k+1}}(0,\bar{i}(\psi^0,\theta
^0)),
\]
it follows that
\[
\frac{1}{\sqrt{n}}\hat{\ell}_{\phi}(\psi^0,\theta^0)\stackrel{\mathcal
{D}}{\rightarrow}N_{d+d_1+d_2+\cdots+d_{k+1}}(0,\bar{i}(\psi^0,\theta
^0)) .
\]

In a similar way, we easily obtain that
\[
-\frac{1}{n}\hat{\ell}_{\phi\phi}(\psi^0,\theta^0)\stackrel{\mathcal
{D}}{\rightarrow}\bar{i}(\psi^0,\theta_j^0).
\]
Therefore, we have that
\[
\sqrt{n}(\hat{\phi}-\phi^0)\stackrel{\mathcal{D}}{\rightarrow
}N_{d+d_1+d_2+\cdots+d_{k+1}}(0,\bar{i}(\psi^0,\theta^0_j)^{-1}),
\]
proving the result.
\end{pf*}

\section*{Acknowledgements}
H. He thanks Professor Peter Hall for his support of this research. The
research of H. He was financially supported by a MASCOS grant from the
Australian Research Council. The research of T.A. Severini was
supported by the U.S. National Science Foundation.

\printhistory

\end{document}